\documentclass[12pt]{amsart}
\usepackage{graphicx}
\usepackage{amsmath,amscd,amssymb}
\usepackage{amsfonts}
\usepackage[dvips]{epsfig}
\usepackage{subfigure}
\usepackage[usenames]{color}
\newtheorem{theorem}{Theorem}[section]

\newtheorem{corollary}[theorem]{Corollary}
\newtheorem{conjecture}[theorem]{Conjecture}
\newtheorem{lemma}[theorem]{Lemma}

\def\la{\lambda}

\def\Om{\Omega}

\newcommand{\Ess}{\textup{ess sup}}

\newcommand{\rz}{\mathbb{R}}

\title{On Riesz Means of Eigenvalues}

\author{Evans M. Harrell II}
\address{
School of Mathematics,
Georgia Institute of Technology, 
Atlanta, GA 30332-0161}
\urladdr{http://www.math.gatech.edu/$\sim$harrell}
\email{harrell@math.gatech.edu}
\author{Lotfi Hermi}
\address{Department of Mathematics,  
University of Arizona, 
Tucson, AZ \, 85721} 
\urladdr{http://math.arizona.edu/$\sim$hermi}
\email{hermi@math.arizona.edu}

\thanks{Manuscript received December 25, 2007; revised XXX XX, XXXX}

\begin{document}

\begin{abstract}
In this article we prove the equivalence of certain inequalities
for Riesz means of eigenvalues of the Dirichlet Laplacian with a
classical inequality of Kac. Connections are made via integral
transforms including those of Laplace, Legendre, Weyl, and Mellin, and the
Riemann-Liouville fractional transform. We also prove new universal 
eigenvalue inequalities and monotonicity principles for 
Dirichlet Laplacians as well as certain
Schr\"odinger operators. 
At the heart of these inequalities are calculations of
commutators of operators, sum rules, and
monotonic properties of Riesz means. In the course of developing
these inequalities we prove new bounds for the 
partition function and the 
spectral
zeta function (cf. Corollaries
3.5-3.7)
and conjecture about additional bounds.
\end{abstract}

%

\vskip 1.0cm

\maketitle


\newpage

\tableofcontents

\section{Riesz Means, Counting Functions, and All That} \label{results}

In \cite{HarHer2} commutator identities
introduced in \cite{HS} were used to
derive both universal and domain-dependent inequalities for
eigenvalues of the Dirichlet Laplacian and the Schr\"odinger operators
with discrete spectra.  
(See also \cite{LevPar}, \cite{AH2},
\cite{AH3}, \cite{SoufiHaIlias}, \cite{HarrellCPDE}.)
In the present article we put those 
notions together with some transform techniques in 
order to connect together several inequalities for spectra, which
have been derived by independent methods in the past.
The essential point is that these inequalities are often 
largely equivalent
under the application of some integral transforms.
Along the way we obtain some improvements and conjecture 
about yet more inequalities.

For the most part we shall concentrate 
on the Dirichlet Laplacian, i.e., on 
the fixed membrane problem on a
bounded domain $\Omega \subset \mathbb{R}^d$,
\begin{equation}
 \Delta u + \lambda \ u =0 \, \text{ in } \Omega,
\end{equation}
subject to Dirichlet boundary condition \[u|_{\partial
\Omega}=0.\]
The boundedness of $\Omega$ serves only to guarantee that the spectrum is purely discrete \cite{Dav2}.
We sometimes treat the Schr\"odinger operator, 
\begin{equation}
- \Delta u + V(x) u = \lambda \ u \, \text{ in } \Omega,
\end{equation}
under circumstances where its spectrum is discrete and bounded from below.
We note that the Schr\"odinger operator $H = - \Delta + V(x) $ may have discrete 
spectrum even when $\Omega$ is not bounded, if $V(x) \rightarrow \infty$ at infinity.  

Eigenvalues are counted with multiplicities and increasingly ordered:
\begin{equation}
\lambda_1 <\lambda_2 \le \ldots \le \lambda_k \le \ldots \to
\infty, \label{ev}
\end{equation}
The eigenvectors, known to form a complete
orthonormal family of $L^2(\Omega)$, are denoted by $u_1, u_2,
\ldots, u_k, \ldots$.

A central object is the Riesz mean of order $\rho>0$. It is
defined, for $z\ge 0$, by
\[R_{\rho}(z) = \sum_{k} \left(z - \lambda_k\right)_{+}^{\rho},\]
where $\left(z-\lambda\right)_{+}:=\max\left(0,z-\lambda\right)$
is the {\em ramp function}.

Here we collect some known properties of $R_{\rho} (z)$ and some
consequences. When $\rho \to 0+$, the Riesz mean reduces to the
counting function (also called the staircase function by
physicists)
\[N(z) = \sum_{\lambda_k \le z} \ 1 = \sup_{\lambda_k \le z} \ k.\]
By convention, this is sometimes written as
\[N(z) = R_{0}(z) = \sum_{k} \left(z - \lambda_k\right)_{+}^0\]
to parallel the definition of the Riesz mean of order $\rho$.
In fact the two are related by the formula
\[R_{\rho}(z) = \int_{0}^{\infty}
\left(z-t\right)_{+}^{\rho}  dN(t) = \rho \int_{0}^{\infty}
\left(z-t\right)_{+}^{\rho-1} N(t) dt.\]

A basic property for $\rho, \delta>0$, sometimes referred
to as {\em Riesz iteration} or as the Aizenman-Lieb procedure
\cite{AL}, is that
\begin{equation}
R_{\rho+ \delta }(\lambda)= \dfrac{\Gamma(\rho+ \delta+1)}
{\Gamma(\rho+1) \ \Gamma(\delta)} \int_{0}^{\infty}
\left(\lambda-t\right)_{+}^{\delta-1} R_{\rho}(t) dt.
\label{riesz1}
\end{equation}

The proof of (\ref{riesz1}) hinges on the Fubini-Tonelli theorem (see p. 3 of 
\cite{ChandraMinak} or \cite{hunder1}) and the fact that
\[\int_{0}^{\infty} \left(1-t\right)_{+}^{p-1} t^{q-1} dt = \dfrac{\Gamma(p) \Gamma(q)}
{\Gamma(p+q)}=\mathcal{B}(p,q),\] where $\Gamma$ and $\mathcal{B}$
denote the Euler functions. Generalizations and further 
facts about Riesz
means can be found in \cite{ChandraMinak} and in some works
related to the Lieb-Thirring inequality (e.g., \cite{HelRob1}
\cite{HelRob2} \cite{hunder1} \cite{hunder2}). We observe that
Riesz iteration is nothing but a Riemann-Liouville fractional
integral transform, the properties of which are tabulated in
\cite{bateman}.

Estimates for these functions of the spectrum have been of
interest for almost a century, since the semiclassical asymptotic formula of
Weyl \cite{Weyl} \cite{AB2} \cite{BH} 
\cite{hunder1} \cite{Kac} \cite{Kac3} \cite{Protter}
for the eigenvalues of the Laplacian,
\begin{equation}
N(z) \sim \dfrac{C_d |\Om| \la^{d/2}}{(2 \pi)^d} = L_{0,d}^{cl} \
|\Omega| z^{d/2} \label{Weyl-asymptotic}
\end{equation}
as $z \to \infty$. Here \begin{equation} L_{0,d}^{cl}:=C_d/(2
\pi)^d \label{classical-constant-0}
\end{equation}
is called the classical constant and $C_d$ is the volume of the
$d$-ball, \[C_d= \pi^{d/2}/\Gamma(1+d/2).\] Note that the Riesz
iteration of (\ref{Weyl-asymptotic}) immediately gives the
statement that
\begin{equation}
R_{\rho}(z) \sim L_{\rho,d}^{cl} \ |\Omega| \ z^{\rho+d/2}
\quad \text{ as } z\to \infty, \label{riesz-weyl}
\end{equation}
where the classical constant is given by
\begin{equation}
L_{\rho,d}^{cl}=\dfrac{\Gamma(1+ \rho)}{\left(4
\pi\right)^{d/2} \Gamma(1+ \rho+ d/2)}. \label{classical}
\end{equation}

\noindent
Furthermore,
\begin{theorem}[Laptev-Weidl] \label{laptev-weidl-thm}
For $\rho \ge 1$, the Riesz means for the Dirichlet Laplacian satisfy
\begin{equation} \label{laptev-weidl} R_{\rho} (z) \le L_{\rho,d}^{cl} \
|\Omega| \ z^{\rho+d/2}.
\end{equation}
\end{theorem}

\par\noindent
\textbf{Remark.} 
In \cite{LapWeidl1} (see also \cite{Lap1} \cite{LapWeidl2})
Laptev and Weidl refer to this as the Berezin-Li-Yau inequality. Indeed, in 1972 Berezin \cite{berezin} proved a general version from which a 1983 inequality of Li-Yau \cite{LY} follows as a corollary (see also \cite{Weidl}).
In terms of the counting function, the Berezin-Li-Yau inequality,
\begin{equation} \label{li-yau-sum}
\sum_{j=1}^k \lambda_j \ge \dfrac{d}{d+2} \dfrac{4 \pi^2
k^{1+2/d}} {\left(C_d |\Om|\right)^{2/d}},
\end{equation}
reads
\begin{equation} \label{li-yau-counting}
N(z) \le \left(\frac{d+2}{d}\right)^{d/2} \, L_{0,d}^{cl} \ |\Om|
z^{d/2}.
\end{equation}
Berezin's version \cite{Safarov} \cite{Lap1} reads
\begin{equation}
\int_0^{z} N(\mu) d \mu \le \dfrac{1}{1+\frac{d}{2}} L_{0,d}^{cl}
z^{1+ d/2} |\Omega|. \label{ber}
\end{equation}
This is just the statement (\ref{laptev-weidl}) for $\rho=1$, 
recalling that the left side is $R_1(z$) and that by \eqref{classical},
\begin{equation} \label{classical-useful}
L_{1,d}^{cl}= \frac{1}{1+\frac{d}{2}} \, L_{0,d}^{cl}.
\end{equation}
Since $N(z)$ is a nondecreasing function, for $\theta>0$,
\begin{equation} \label{berez}
N(z) \le \frac{1}{\theta z} \displaystyle{\int_{z}^{(1 + \theta)
z} N(\mu) d \mu} \le \frac{1}{\theta z} \displaystyle{\int_0^{(1 +
\theta) z} N(\mu) d \mu} \le \dfrac{\left(1+
\theta\right)^{1+d/2}}{\left(1+\frac{d}{2}\right) \theta}
L_{0,d}^{cl} |\Omega| z^{d/2}. \notag
\end{equation}
The Berezin-Li-Yau bound (\ref{li-yau-counting}) follows by
setting $\theta=2/d$. In a rather straightforward way, the method 
of \cite{Lap1} and \cite{Safarov} for proving \eqref{li-yau-counting}
yields a formula that 
interpolates between Berezin-Li-Yau
($\rho=0$) and Laptev-Weidl ($\rho\ge 1$). 
\begin{theorem} \label{obvious} For $0 \le \rho<1$,
the Riesz means for the Dirichlet Laplacian satisfy
\begin{equation}  \label{obvious-result}
R_{\rho} (z) \le K_{\rho, d} \, \Gamma(1+ \rho) \,
\Gamma(2-\rho) L_{1, d}^{cl} |\Omega| \, z^{\rho+\frac{d}{2}}
\end{equation}
where \[K_{\rho, d}= \inf_{\theta>0}
\dfrac{\left(1+\theta\right)^{1+d/2}}{\theta^{1-\rho}}=
\dfrac{1}{\left(1-\rho\right)^{1-\rho}} \,
\dfrac{\left(1+d/2\right)^{1+d/2}}{\left(\rho+d/2\right)^{\rho+d/2}}.\]
\end{theorem}
\par\noindent \textbf{Proof.} For the range of values of $\rho$ considered,
$R_{\rho} (z)$ is a nondecreasing function of $z>0$. Therefore,
for $\theta>0$ and $\delta>0$, 
\begin{eqnarray}
\left(\theta z\right)^{\delta} R_{\rho}(z) &\le & \delta \, 
\displaystyle{\int_{z}^{(1 + \theta) z} \left(z+\theta z-
t\right)^{\delta-1} R_{\rho}(t) dt} \\ \notag &\le & \delta \,
\displaystyle{\int_{0}^{(1 + \theta) z} \left(z+\theta z-
t\right)^{\delta-1} R_{\rho}(t) dt}
\\ \notag
&=&
\dfrac{\Gamma(\rho+1)
\Gamma(\delta+1)}{\Gamma(\rho+\delta+1)} \, R_{\rho+\delta} (z
+ \theta z). \notag
\end{eqnarray}

Therefore
\begin{equation}
R_{\rho}(z) \le \dfrac{\Gamma(\rho+1) \Gamma(\delta+1)}
{\Gamma(\rho+\delta+1)} \, \inf_{\theta>0}
\dfrac{R_{\rho+\delta} \left((1+\theta)z)\right)}{\left(\theta
z\right)^{\delta}}.
\label{useful-obvious}
\end{equation}
Specializing to the case $\rho+\delta=1$ and using Berezin-Li-Yau (\ref{laptev-weidl}) leads to
\eqref{obvious-result}.
The optimal bound occurs when
$\theta=\dfrac{1-\rho}{\rho+d/2}$. This reduces to the
estimate \eqref{li-yau-counting} when $\rho \to 0+$. \hfil \qed

\begin{figure}[htb!]
  \begin{center}
\caption{\textit{Comparison of the constant in \eqref{obvious-result} and the Riesz iteration of \eqref{li-yau-counting} with the classical constant $L_{\rho, d}^{cl}$ for $0\le \rho \le 1$ and $d=3$.}}
    \label{figBLY}
    \includegraphics[width=.95 \textwidth, bb=0 0 320 197 ]{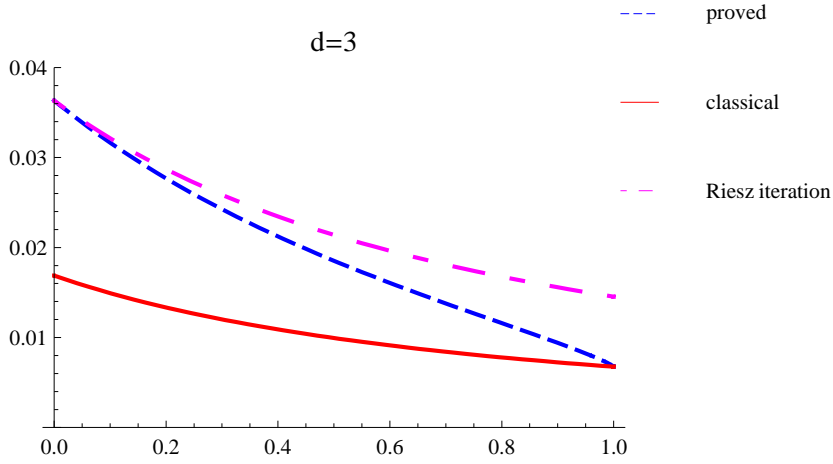}
  \end{center}
\end{figure}

Fig.~ \ref{figBLY} depicts the interpolation between
Berezin-Li-Yau ($\rho=0$) and Laptev-Weidl ($\rho\ge 1$) in dimension 3, as well as the resulting graph from direct Riesz iteration of \eqref{li-yau-counting} which results in a weaker bound (see also the discussion in \cite{Weidl}). 

By developing ideas from 
\cite{HS}, \cite{AH3} it was proved in \cite{HarHer2}
that for $\rho\ge 2$,
\begin{equation} \label{ha-stu2}
\sum_k \left(z-\lambda_k\right)_{+}^{\rho} \le \dfrac{2
\rho}{d} \, \sum_k \lambda_k \,
\left(z-\lambda_k\right)_{+}^{\rho-1},
\end{equation}
Thereby extending the 
``Yang-type inequality'' \cite{Yang}, \cite{A1} (see also \cite{HS} \cite{LevPar} \cite{AH2} and the appendix to \cite{ChengYang2}), {\it viz.,}
\begin{equation} \sum_k
\left(z-\lambda_k\right)_{+}^{2} \le \dfrac{4}{d} \, \sum_k
\lambda_k \, \left(z-\lambda_k\right)_{+}, \label{yang}
\end{equation} 
corresponding to $\rho=2$.
In Section \ref{harrell-stubbe} we shall show how the
inequalities for $\rho > 2$ can be directly deduced from
\eqref{yang}.

Two final functions of the spectrum will be of interest, 
the spectral zeta function defined by
\[\zeta_{spec} (\rho) =\sum_{k=1}^{\infty} \dfrac{1}{\lambda_k^{\rho}},\]
and the partition function (= trace
of the heat kernel) $Z(t)$.  We recall the asymptotic formula 
of Kac \cite{Kac} for $Z(t)$: 

\begin{equation} \label{Kac-asymptotic}
Z(t) = \sum_{k=1}^{\infty} e^{-\lambda_k t}  \sim
\dfrac{|\Omega|}{\left(4 \pi t\right)^{d/2}},
\end{equation}
and observe that it can be proved with an application of the
Laplace transform $\mathcal{L} \{ f \}(t) = \int_{0}^{\infty} f(z)
e^{-z t} d z$ to (\ref{Weyl-asymptotic}).
In \cite{Kac3} Kac also used ``the principle of not feeling the boundary''
to derive the inequality
\begin{equation} \label{Kac-inequality}
Z(t) = \sum_{k=1}^{\infty} e^{-\lambda_k t}  \le
\dfrac{|\Omega|}{\left(4 \pi t\right)^{d/2}}.
\end{equation}
In \cite{HS} this was improved to the statement that 
$t^{d/2} Z(t)$ is a nonincreasing function that saturates when $t \to 0+$.

We remark here on some extensions to the case of Schr\"odinger operators.
If the potential function $V(x) \ne 0$, then eq. \eqref{yang} becomes
\begin{equation} 
\sum_k
\left(z-\lambda_k\right)_{+}^{2} \le \dfrac{4}{d} \, \sum_k
T_k \, \left(z-\lambda_k\right)_{+}, \label{Schryang}
\end{equation} 
where 
\begin{equation}
T_k := \int_{\Omega}{|\nabla u_k|^2} = \lambda_k - \int_{\Omega}{V(x) |u_k|^2} := 
\lambda_k - V_k,
\label{Tk}
\end{equation} 
(cf. eq. (12) of \cite{HS}).  As was remarked in \cite{HS}, $T_k$ is often bounded above by a multiple of $\lambda_k$ under general assumptions on $V$, for
example those that guarantee a virial inequality.  Another circumstance in which such a
bound is possible is when the negative part of $V$ is relatively bounded by the Laplacian \cite{Kato}, \cite{ReSi}, \cite{Dav2}, whether in the sense of operators or of quadratic forms.  As an example, according to the Gagliardo-Nirenberg inequality (e.g., \cite{Allegretto}), there is a dimension-dependent constant $K_{GN,d}$, such that if  
\par\noindent
$V_- := \max(0, -V(x)) \in L^{d/2}$, then
$\int_{\Omega}{V_-(x) |u_k|^2} \le K_{GN,d} \, \|V_-\|_{d/2} \, T_k.$
Under these circumstances, \[T_k \le \lambda_k 
+ \int_{\Omega}{V_-(x) |u_k|^2} \le \lambda_k + K_{GN,d} \|V_-\|_{d/2}  \, T_k.\]
If, moreover, $ \|V_-\|_{d/2} < 1/K_{GN,d}$, then it follows that
\[T_k < \frac{1}{1 - K_{GN,d} \|V_-\|_{d/2}} \, \lambda_k.\] In \cite{Allegretto}, the constant $K_{GN,d}$ is given in the explicit form
\[K_{GN,d}=\frac{(d-1)^2}{(d-2)^2 \, d},\]
thus restricting the dimension to $d \ge 3$.

Because there are many circumstances where a bound of this form applies, for future purposes we refer to:
\medskip

\noindent
\textbf{Assumption $\Sigma$.}  For some $\sigma < \infty$,
$T_k \le \sigma \, \lambda_k$.

\medskip
The article is organized as follows. We first prove the equivalence
of several old and new inequalities for the spectrum of the
Dirichlet Laplacian. Central to our argument is a monotonicity
principle proved in \cite{HarHer2}, to which we offer a new path
via integral transforms. We then use a sum rule in the style of Bethe 
\cite{BetJac} \cite{Jac} to recover
bounds which compete with the Berezin-Li-Yau inequality 
\eqref{laptev-weidl}, and also with results recently proved in \cite{HarHer2}.
Finally we comment on some possible
corrections to the Berezin-Li-Yau inequality and related inequalities. 

\section{The Equivalence of Several Inequalities for spectra}

In this section we show that many universal and geometric bounds
for spectra of the Dirichlet Laplacian,
which have been proved in the literature by independent methods, 
may in fact be derived from one another by the application of 
the Laplace transform and some classical 
inequalities.  In particular, for 
$\rho \ge 2$, it will be shown that
the Kac inequality (\ref{Kac-inequality}) and the
Berezin-Li-Yau inequality (\ref{laptev-weidl}) are
equivalent by the Laplace transform.  These inequalities are 
seen to be corollaries of the Riesz-mean inequalities of \cite{HS} \cite{HarHer2}, which in turn can all be derived from the case $\rho=2$, originating with 
Yang. 

With some minor modifications, similar inequalities are then proved 
for Schr\"odinger spectra.

\subsection{Kac from Berezin-Li-Yau}
For the pure Laplacian, with no added potential, we start 
by showing that the Kac inequality (\ref{Kac-inequality}) can be derived from the
Berezin-Li-Yau inequality (\ref{laptev-weidl}) as
an alternative to Kac's
``principle of not feeling the boundary''. 
Begin with the
observation that the Laplace transform yields
\begin{equation}
\mathcal{L} \left( \left(z- \lambda_k\right)_{+}^{\rho}\right)=
\dfrac{\Gamma(\rho+1) \ e^{-\lambda_k \ t}}{t^{\rho+1}}.
\label{useful}
\end{equation}
Applying this to (\ref{laptev-weidl}) immediately leads to
\[\dfrac{\Gamma(\rho+1)}{t^{\rho+1}} Z(t) \le L_{\rho, d}^{cl} \, |\Omega|
\dfrac{\Gamma(\rho+1+\frac{d}{2})}{t^{\rho+1+\frac{d}{2}}},\]
which upon simplification reads
\[Z(t) \le \dfrac{|\Omega|}{t^{\frac{d}{2}}} \,
\dfrac{L_{\rho, d}^{cl} \,
\Gamma(\rho+1+\frac{d}{2})}{\Gamma(\rho+1)}.\] Using the
definition of $L_{\rho, d}^{cl}$ in (\ref{classical}) results in
(\ref{Kac-inequality}).  Indeed it is only necessary to have 
\eqref{laptev-weidl} for a single value of $\rho$.

We observe that the same argument relates the 
Kac-Ray inequality \cite{Kac} \cite{Kac3} \cite{Ray} \cite{vdB}, 
\begin{equation} \label{kac-ray}
Z(t) \le \frac{1}{(4 \pi t)^{d/2}} \,
\displaystyle{\int_{\rz^d} e^{-t V(x)} dx}
\end{equation}
(also known in the literature as the Golden-Thompson inequality \cite{DFLP})
to the Lieb-Thirring inequality \cite{Lap2} \cite{LapWeidl2}
\begin{equation} \label{lt-sharp}
R_{\rho}(z) \le L_{\rho, d}^{cl} \, \displaystyle{\int_{\rz^d} \left(z-V(x)\right)_{+}^{\rho+d/2} dx },
\end{equation}
for the Laplace transform of \eqref{lt-sharp} yields \eqref{kac-ray}.

\subsection{Kac from Yang} \label{laplace-yang} 
Next we show how to obtain Kac's inequality
\eqref{Kac-inequality} directly from Yang's inequality
\eqref{yang} and the asymptotic formula \eqref{Kac-asymptotic}.
The link is a result of Harrell and Stubbe \cite{HS}:

\begin{theorem}
The function $t^{d/2} \ Z(t)$
is a nonincreasing function.
\label{HS-theo}
\end{theorem}

In \cite{HS}, this theorem was derived from a trace
identity, but here we show that it can alternatively be proved from
Yang's inequality
\eqref{yang}. 

Apply the Laplace transform to both sides of
(\ref{yang}), written now in the form
\[\sum_{k=1}^{\infty}
\left(z-\lambda_k\right)_{+}^{2} \le \dfrac{4}{d} \,
\sum_{k=1}^{\infty} \lambda_k \, \left(z-\lambda_k\right)_{+},
\]
and use (\ref{useful}) to obtain the differential inequality
\begin{equation} \label{diff}
Z(t)\le -\dfrac{2}{d} \, t \, Z'(t)
\end{equation}
or, after combining,
\[\left(t^{d/2} \, Z(t) \right)^{'} \le 0.\]
Kac's inequality is then immediate, employing
\eqref{Kac-asymptotic} in the form
\begin{equation}
\lim_{t \to 0+} \ t^{d/2} Z(t) = \dfrac{|\Omega|}{\left(4 \pi
\right)^{d/2}}. \label{well-known}
\end{equation}

\subsection{Riesz-mean inequalities for $\rho > 2$ from Yang}
\label{harrell-stubbe}

In this section we show how to prove
\eqref{ha-stu2}
directly from \eqref{yang}. 

\begin{theorem} \label{ha-stu-thm}
\cite{HS}
For $\rho \ge 2$ and $z \ge 0$,
\begin{equation}
R_{\rho} (z) \le \dfrac{\rho}{\rho+\frac{d}{2}} \, z \
R_{\rho-1} (z). \label{ha-stu1}
\end{equation}
\end{theorem}

As in the original proof from first principles \cite{HarHer2}, we note that 
\eqref{ha-stu1} is equivalent to \eqref{ha-stu2}.
To see this,
rewrite $\lambda_k \left(z-\lambda_k\right)_{+}^{\rho-1}$
in (\ref{ha-stu2}) as
\[\left(-z+\lambda_k +z \right) \,
\left(z-\lambda_k\right)_{+}^{\rho-1},\] 
and rearrange terms.


In order to use Riesz iteration we now rewrite (\ref{yang}) 
for $t \le z$ as
\begin{equation}
\sum_k \left(z - \lambda_k-t\right)_{+}^2 \leq \dfrac{4}{d} \sum_k
\lambda_k \left(z - \lambda_k-t\right)_{+}.\notag
\end{equation}
Multiply both sides by $t^{\rho-3}$, 
and then integrate between $0$ and
$\infty$. By (\ref{useful}), there results
\begin{equation}
\sum_k \left(z - \lambda_k\right)_{+}^{\rho} \leq \dfrac{2}{d} \
\dfrac{\Gamma(\rho+1) \Gamma(2)}{\Gamma(\rho) \Gamma(3)}
\sum_k \lambda_k \left(z - \lambda_k\right)_{+}^{\rho-1}. \notag
\end{equation}
With $\Gamma(\rho+1) = \rho \ \Gamma(\rho)$, this simplifies to
\begin{equation} \label{dfq}
\sum_k \left(z - \lambda_k\right)_{+}^{\rho} \leq
\dfrac{2\rho}{d} \ \sum_k \lambda_k \left(z -
\lambda_k\right)_{+}^{\rho-1}, 
\end{equation}
which is the statement of Theorem \ref{ha-stu-thm}. 
It was shown in 
\cite{HarHer2} that
\eqref{dfq} is equivalent to the differential inequality
\begin{equation} \label{dfq3}
R_{\rho} (z) \le \dfrac{1}{\rho+\frac{d}{2}} \, z \,
{R'}_{\rho}(z),
\end{equation}
and hence to a monotonicity principle,

\begin{theorem}[\cite{HarHer2}] \label{mono}
The function \[z \mapsto
\dfrac{R_{\rho}(z)}{z^{\rho+\frac{d}{2}}}\] is a nondecreasing
function of $z$, for $\rho\ge 2$.
\end{theorem}

\par\noindent \textbf{Remark.} In \cite{AH3}, it was proved that if
$\gamma_m(\rho)$ is the unique solution of \[ \sum_k \left(z -
\lambda_k\right)_{+}^{\rho} = \dfrac{2\rho}{d} \ \sum_k
\lambda_k \left(z - \lambda_k\right)_{+}^{\rho-1}.\] for $z\ge
\lambda_{m}$, then $\lambda_{m+1}\le \gamma_m(\rho)$. Moreover
$\lambda_{m+1}\le \gamma_m(\rho) \le \gamma_m(\rho')$ for
$2\le \rho \le \rho'$.
Given that the cases $\rho > 2$ of (\ref{ha-stu1}) follows from the case 
$\rho=2$, it might be thought that it is not sharp for large $\rho$.  
To the contrary, it was shown in \cite{HarHer2} that (\ref{ha-stu1}) implies 
strict bounds with the correct power corresponding to Weyl's law.  Indeed:

\begin{theorem} \label{ha-stu-thm2}
The constant in inequality (\ref{ha-stu1})  for $\rho\ge 2$
cannot be improved.
\end{theorem}
\par\noindent \textbf{Proof.} The proof proceeds by contradiction. Suppose there
exists a constant $C(\rho, d)<\dfrac{\rho}{\rho +
\frac{d}{2}}$ such that
\begin{equation}
R_{\rho}(z) \leq
 \ C(\rho, d) \ z \ R_{\rho-1}(z) . \label{harrell-stubbe-berezin}
\end{equation}
Dividing both sides by $z^{\rho +\frac{d}{2}} \ |\Omega|$, then
sending $z\to \infty$, leads to
\begin{equation}
L_{\rho, d}^{cl} \leq C(\rho, d) \ L_{\rho-1,d}^{cl}.\notag
\end{equation}
However,
\[L_{\rho, d}^{cl}=\dfrac{\rho}{\rho+ \frac{d}{2}} \ L_{\rho-1,d}^{cl},\]
and therefore $C(\rho, d) \ge \dfrac{\rho}{\rho + \frac{d}{2}}$.
This contradicts the assumption and proves the claim. \hfill\qed


\subsection{Berezin-Li-Yau from Harrell-Stubbe}
At this stage we make the simple observation that 
for $\rho \ge 2$, the Berezin-Li-Yau inequality
(\ref{laptev-weidl}) follows immediately from
inequality (\ref{ha-stu2}) (or (\ref{ha-stu1}))
by virtue of the monotonicity principle of
Theorem \ref{mono} and the asymptotic
formula \eqref{riesz-weyl}.

\subsection{Riesz-mean inequalities for $\rho <
2$ from Yang}\label{har-stub2}

In \cite{HarHer2} the difference inequality
\begin{equation}\label{rho-small}
\sum_{k} \left(z-\lambda_k\right)_{+}^{\rho} \le \frac{4}{d} \,
\sum_{k} \lambda_k \, \left(z-\lambda_k\right)_{+}^{\rho-1}
\end{equation}
for $1 < \rho \le 2$ was obtained from first principles and used to prove Weyl-type
universal bounds for ratios of eigenvalues.
Eq. \eqref{rho-small} implies a differential inequality and monotonicity
principle similar to Theorem \ref{mono}, but as an alternative we show how to obtain \eqref{rho-small} using the ``Weighted Reverse Chebyshev Inequality''  (see, for example, p.~43 of \cite{HLP} or \cite{AH3}):

\begin{lemma} Let $\{a_i\}$ and $\{b_i\}$ be two real sequences,
one of which is nondecreasing and the other nonincreasing,
and let $\{w_i\}$ be a sequence of nonnegative weights. Then, 
\begin{align}
\sum_{i=1}^m w_i \ \sum_{i=1}^m w_i \ a_i b_i \le \sum_{i=1}^m w_i
\ a_i \ \sum_{i=1}^m w_i \ b_i. \label{eq:4.4}
\end{align}
\end{lemma}

Making the choices
$w_i= \left(z - \lambda_k\right)_{+}^{\rho_1}$, $a_i
= \frac{\lambda_k}{\left(z - \lambda_k\right)_{+}}$, and
$b_i=\left(z - \lambda_k\right)_{+}^{\rho_2-\rho_1}$
with $\rho_1\le \rho_2 \le 2$, the conditions of the 
lemma are satisfied and we get 

\begin{align}
\sum_{k} \left(z - \lambda_k\right)_{+}^{\rho_1} \ \sum_{k}
\left(z - \lambda_k\right)_{+}^{\rho_2-1} \ \lambda_k \le
\sum_{k} \left(z - \lambda_k\right)_{+}^{\rho_2} \ \sum_{k}
\left(z - \lambda_k\right)_{+}^{\rho_1-1} \ \lambda_k, \notag
\end{align}
which is equivalent to
\begin{align}
\dfrac{\sum_{k} \left(z -
\lambda_k\right)_{+}^{\rho_1}}{\sum_{k} \left(z -
\lambda_k\right)_{+}^{\rho_1-1} \ \lambda_k} \le \dfrac{\sum_{k}
\left(z - \lambda_k\right)_{+}^{\rho_2}} {\sum_{k} \left(z -
\lambda_k\right)_{+}^{\rho_2-1} \ \lambda_k}.  \label{ha-stub2-1}
\end{align}

To obtain inequality \eqref{rho-small}, now set
$\rho_1=\rho$ and $\rho_2=2$ in the above and use
inequality \eqref{yang} to 
estimate the right side. 
We observe that
inequality \eqref{rho-small} not only 
implies familiar results for
$\rho=1$ and $\rho=0$ (the  Hile-Protter inequality \cite{HP}),
but also hitherto unexplored inequalities for $\rho<0$. 

\subsection{Berezin-Li-Yau from Kac, for $\rho \ge 2$} \label{hs2bly}

We showed earlier how to obtain Kac's inequality
\eqref{Kac-inequality} from
\eqref{laptev-weidl}. In this section, we show the reverse, and thus the full equivalence of the two statements.
Throughout this section we assume
$\rho \ge 2$. 

As a result of the Monotonicity Theorem \ref{mono},
for $z\ge z_0$,
\begin{equation} \label{mono1}
R_{\rho} (z) \ge R_{\rho} (z_0) \,
\left(\frac{z}{z_0}\right)^{\rho+ d/2}.
\end{equation}
With $\mu = -z_0+z > 0$,
\begin{equation} \label{mono2}
R_{\rho} (\mu +z_0) \ge R_{\rho} (z_0) \,
\left(\frac{\mu+z_0}{z_0}\right)^{\rho+ d/2}.
\end{equation}
The Laplace transform  of a shifted function is given by
the formula (see p.~3 of \cite{RobKauf})
\[
\mathcal{L} \left( f(\mu+ z_0)\right) = e^{z_0 \, t}
\left(\mathcal{L}(f) - \int_{0}^{z_0} e^{-t \mu} f(\mu) d\mu
\right)
\]
We apply the Laplace transform to (\ref{mono2}),
noting that for the left side,
\begin{equation} \label{mono3}
\mathcal{L} \left( \left(\mu+ z_0-\lambda_k\right)_{+}^{\rho}
\right) =e^{\left(z_0-\lambda_k\right)_{+} t}
\left(\frac{\Gamma(\rho + 1)}{t^{\rho+1}} - \
\int_{0}^{\left(z_0-\lambda_k\right)_{+} t} e^{-t \mu}
\mu^{\rho} d\mu \right),
\end{equation}
whereas on the right,
\begin{equation} \label{mono4}
\mathcal{L} \left( \left(\mu+ z_0\right)^{\rho+d/2} \right)
=e^{z_0 \, t} \left(\frac{\Gamma(\rho + 1+ d/2)}{t^{\rho+ 1+
d/2}} - \int_{0}^{z_0 \, t} e^{-t \mu} \mu^{\rho+d/2} d\mu
\right).
\end{equation}
We note the appearance of the incomplete Gamma function (see p. 260 of \cite{AS})
\[\gamma(a,x) = \int_{0}^{x} e^{-\mu} \mu^{a-1} d \mu.\]
Putting these facts together, we are led to
\begin{eqnarray}
\sum_{k} e^{\left(z_0-\lambda_k\right)_{+} t} \, \Big\{
\frac{\Gamma(\rho+1)}{t^{\rho+1}} & -
&\gamma\left(\rho+1,(z_0-\lambda_k)_{+} t \right) \Big\} \ge
\notag \\ & &\frac{R_{\rho}(z_0)}{z_0^{\rho+d/2}} e^{z_0 \, t}
\, \Big\{\frac{\Gamma(\rho+1+d/2)}{t^{\rho+1+d/2}} -
\gamma(\rho+1+d/2, z_0 \, t) \Big\}. \notag
\end{eqnarray}
We now notice that
\begin{equation} \label{mono-expo}
\sum_{k} e^{\left(z_0-\lambda_k\right)_{+} t} \le e^{z_0 \, t} \
\sum_{k=1}^{\infty} e^{-\lambda_k t}  = e^{z_0 \, t} \, Z(t).
\end{equation}
Therefore, after a little simplification,
\begin{equation} \label{mono-main}
\frac{\Gamma(\rho+1)}{\Gamma(\rho+1+d/2)} \, t^{d/2} Z(t) \ge
\frac{R_{\rho}(z_0)}{z_0^{\rho+d/2}} + \mathcal{R}(t),
\end{equation}
where the remainder term $\mathcal{R}(t)$ has the explicit form
\begin{eqnarray*}
\mathcal{R}(t)= \frac{t^{d/2}}{\Gamma(\rho+1+d/2)} & e^{-z_0 t}
& \sum_{k} e^{\left(z_0-\lambda_k\right)_{+} t} \,
\gamma(\rho+1,\left(z_0-\lambda_k\right)_{+} t )
\\ \notag
&-& \frac{t^{d/2}}{\Gamma(\rho+1+d/2)} \,
\frac{R_{\rho}(z_0)}{z_0^{\rho+d/2}} \, \gamma(\rho+1+d/2,
z_0 \, t)
\end{eqnarray*}
Notice that $\lim_{t\to 0} \mathcal{R}(t)=0$. Sending $t\to 0$ in (\ref{mono-main}) and again incorporating
(\ref{well-known}) leads to
\begin{equation} \label{mono-main2}
\frac{\Gamma(\rho+1)}{\left(4 \pi\right)^{d/2} \,
\Gamma(\rho+1+d/2)} \, |\Omega| \ge
\frac{R_{\rho}(z_0)}{z_0^{\rho+d/2}} .
\end{equation}
We finish by observing that the constant on the left side of
\eqref{mono-main2} is the
classical constant $L_{\rho,d}^{cl}$ from (\ref{classical}).
Hence Berezin-Li-Yau follows for $\rho\ge 2$, as claimed.  
In summary,
when $\rho \ge 2$ the Berezin-Li-Yau inequality is equivalent
to the Kac inequality.

\subsection{Extension to Schr\"odinger spectra}

We have shown above that a family of 
universal inequalities and monotonicity theorems
for Riesz means of Laplace spectra
can be derived from \eqref{yang}.  Under
Assumption $\Sigma$, 
{\it viz}., $T_k \le \sigma \, \lambda_k$, for a constant 
$\sigma<\infty$, 
a similar inequality, differing 
from \eqref{yang} only by the value of a constant,
holds for Schr\"odinger operators. 
Consequently, the 
universal inequalities and monotonicity theorems
discussed above continue to hold for
$H = -\Delta + V(x)$,
with appropriately adjusted constants.  

\begin{theorem}
Assume that $H = - \Delta + V(x)$ is essentially self-adjoint on 
$C_c(\Omega)$; has purely discrete spectrum with $\lambda_1 > -\infty$;
and satisfies Assumption $\Sigma$.  Then

a) (Riesz means, $\rho \ge 2$)   For  $\rho \ge 2$ and $z \ge 0$,
\begin{equation}\label{mod-sigma-large}
R_{\rho}(z) \le \frac{\rho}{\rho + \frac{d}{2 \sigma}} z R_{\rho-1}(z),
\end{equation}

and consequently
the function \[z \mapsto
\dfrac{R_{\rho}(z)}{z^{\rho+\frac{d}{2 \sigma}}}\] 
is a nondecreasing
function of $z$.

b) (Riesz means, $\rho \le 2$)   For  $1 < \rho \le 2$ and $z \ge 0$,
\begin{equation}\label{mod-sigma-small}
R_{\rho}(z) \le \frac{1}{1 + \frac{d}{4 \sigma}} z R_{\rho-1}(z),
\end{equation}

and consequently
the function \[z \mapsto
\dfrac{R_{\rho}(z)}{z^{\rho+\frac{d \rho}{4 \sigma}}}\] 
is a nondecreasing
function of $z$.
\end{theorem}

\par\noindent \textbf{Remarks.}

1.
The proofs are precisely like the ones given above with a change of constant, and are therefore omitted.  The assumptions in the theorem suffice to allow 
eq. (12) of \cite{HS} as a replacement for \eqref{yang} (see also \cite{HarHer2}).

2.
With a similar argument, a modification of Kac's
inequality was obtained in 
\cite{HS}:  {\it The function $t^{d/2 \sigma} Z(t)$ is monotonically
nonincreasing in $t$.}

3.
We recall the values of $\sigma$ in three 
simple situations in which Assumption $\Sigma$ holds:
\begin{itemize}
\item[(i)]  If $V(x) \ge 0$, then $\sigma = 1$. 
\item[(ii)]  If for some $\beta > 0$, ${\bf x} \cdot \nabla V(x) \le \beta \, V(x)$, 
then $\sigma = \beta/(2+\beta)$  (cf. \cite{HS}).
\item[(iii)]  As in the earlier discussion, if $\|V\|_{d/2} < K_{GN,d}$, then $\sigma = \frac{1}{1-\|V\|_{d/2}/ K_{GN,d}}$.
\end{itemize}

\section{Lower Bounds for Riesz Means, Zeta Functions, and Partition Functions}
\label{lower-bounds}

In this section, we obtain lower bounds on $R_{\rho}(z)$, which for some
parameter values improve the lower bounds obtained in \cite{HarHer2}.  As corollaries we get lower bounds on spectral zeta functions and on the partition function.

\begin{theorem} \label{riesz-low-thm}
For $\rho\ge 1$
\begin{equation}
R_{\rho}(z) \ge H_d^{-1}  \
\dfrac{\Gamma(1+\rho) \Gamma(1+d/2)}{\Gamma(1+ \rho+ d/2)} \
\lambda_1^{-d/2} \left(z-\lambda_1\right)_{+}^{\rho+d/2}. \label{riesz-low}
\end{equation}
\end{theorem}
Here
\begin{equation}H_d=\dfrac{2 \ d}{j_{d/2-1,1}^2
J_{d/2}^2(j_{d/2-1,1})} \label{constant-riesz-low}
\end{equation}
is a universal constant which depends on the dimension $d$, while
$J_n(x)$ and $j_{n,p}$ denote, respectively, the Bessel function
of order $n$, and the $p$th zero of this function (see \cite{AS}).
The case $\rho=1$ of (\ref{riesz-low}) has been proved in
\cite{Her} using the Rayleigh-Ritz method, 
and in \cite{Safarov} Safarov derived similar lower bounds, with a lower constant.
Yet another independent proof and generalization appeared in \cite{FLM}, in the spirit of \cite{Lap1}. We shall obtain some improvement by use of Riesz iteration and Chiti's isoperimetric lemma \cite{Chi3}. Note that ineq. \eqref{riesz-low} is valid for both the eigenvalues of the Dirichlet-Laplacian and the class of Schr\"odinger operators treated in this article (cf. \cite{FLM}).

The starting point is the {\it Bethe sum rule} as it appears in \cite{LevPar}:

\begin{equation} \label{bethe}
\sum_{k} \left(\lambda_k - \lambda_j\right) |a_{jk}(\xi)|^2 = |\xi|^2,
\end{equation}
where
\begin{equation}
a_{jk}(\xi) =\int_{\Omega} u_k u_j e^{i x \cdot \xi} dx,
\label{bethe-coeff}
\end{equation}
and $\xi \in \rz^d$. 

The Bethe sum rule provides an elementary proof of a lemma of Laptev \cite{Lap1}, originally proved using pseudodifferential calculus:

\begin{theorem}[Laptev \cite{Lap1}] \label{low-bound-laptev}
\begin{equation}
\label{laptev-2} \sum_{j} \left(z- \lambda_j \right)_{+} \ge L_{1,
d}^{cl} \, {\tilde{u}_{1}}^{-2} \,
\left(z-\lambda_1\right)_{+}^{1+d/2}.
\end{equation}
where $\tilde{u}_{1}=\textup{ess sup} |u_1|$ and $L_{1, d}^{cl}$
is given in (\ref{classical-constant-0}).
\end{theorem}

\par\noindent \textbf{Remarks.}
Laptev's form of the inequality reads
\begin{equation}
\label{laptev-1} \sum_{j} \left(z- \lambda_j \right)_{+} \ge
\frac{1}{1+\frac{d}{2}} \, L_{0, d}^{cl} \, {\tilde{u}_{1}}^{-2}
\left(z-\lambda_1\right)_{+}^{1+d/2},
\end{equation}
which is equivalent by dint of 
(\ref{classical-useful}).

\par\noindent \textbf{Proof.} In (\ref{bethe}), choose $j=1$, to get
\[
\sum_{k} \left(\lambda_k - \lambda_1\right) |a_{1k}(\xi)|^2 = |\xi|^2.
\]
Let $z>\lambda_1$.  One can always find an integer 
$N$ such that
\[\lambda_N < z \le \lambda_{N+1},\]
allowing the sum to be split as
$\sum_k =\sum_{k=1}^N + \sum_{k=N+1}^{\infty}$.
We can replace each term in
$\sum_{k=N+1}^{\infty}{(\dots)}$ by \[\left(z- \lambda_1\right) \
|a_{1k}(\xi)|^2.\] Hence
\begin{equation} \label{step1-low}
\sum_{k=1}^N
\left(\lambda_{k}- \lambda_1\right) \ |a_{1k}(\xi)|^2 + \left(z-
\lambda_1\right) \ \left(1-\sum_{k=1}^N |a_{1k}(\xi)|^2 \right) \le
|\xi|^2.
\end{equation}
Here we have exploited the completeness of the orthonormal family $\{u_k\}_{k=1}^{\infty}$, 
noting that
\[\sum_{k=1}^{\infty} |a_{1k}(\xi)|^2 = \int_{\Omega} |u_1 e^{i x \cdot \xi}|^2=1.\]
Therefore
\[\sum_{k=N+1}^{\infty} |a_{1k}(\xi)|^2 = 1- \sum_{k=1}^N |a_{1k}(\xi)|^2.\]
These identities reduce (\ref{step1-low}) to
\begin{equation}
\left(z - \lambda_1\right)_{+} \le |\xi|^2 + \sum_{k} \left(z -
\lambda_k \right)_{+} |a_{1k}(\xi)|^2. \label{step2-low}
\end{equation}
(The statement is true by default for $z \le \lambda_1$.)
One then integrates over a ball $B_r \subset \rz^d$ of radius $r$. To
simplify the notation we use
\[|B_r| =\text{volume of } B_r= C_d \, r^d,\]
and
\[I_2(B_r) = \int_{B_r} |\xi|^2 d \xi= \frac{d}{d+2} C_d \, r^{d+2}.\]
Ineq. (\ref{step2-low}) reduces to
\begin{equation}
\left(z - \lambda_1\right)_{+} \le \dfrac{I_2(B_r)}{|B_r|} +
\sum_{k} \left(z - \lambda_k \right)_{+} \dfrac{\int_{B_r}
|a_{1k}(\xi)|^2 d\xi}{|B_r|}. \label{step3-low}
\end{equation}
By the Plancherel-Parseval identity
\begin{eqnarray}
\dfrac{1}{(2 \pi)^{d}} \ \int_{B_r} |a_{1k}(\xi)|^2  d\xi  & \le
& \int_{\Om} |u_1|^2 |u_k|^2 dx \\ \notag & \le & \textup{ess sup}
|u_1|^2 \int_{\Om}  |u_k(x)|^2 dx
\\ \notag
& =  & \textup{ess sup} |u_1|^2.\\
\label{parseval}
\end{eqnarray}
Incorporating (\ref{parseval}) into (\ref{step3-low}) and simplifying the expression
leads to
\begin{equation}
\sum_{k} \left(z - \lambda_k \right)_{+} \ge \tilde{u}_{1}^{-2}
L_{0,d}^{cl} \, r^d \, \left[ \left(z - \lambda_1\right)_{+} -
\frac{d}{d+2} r^2 \right]. \label{step4-low}
\end{equation}
Optimizing over $r$ results in the statement of the theorem. \hfill\qed

As an immediate consequence of Theorem \ref{low-bound-laptev} and Riesz
iteration, we have the following.
\begin{corollary} \label{cor-univer} For $\rho \ge 1$
\begin{equation}
\label{riesz-lap} \sum_{k} \left(z-\lambda_k\right)_{+}^{\rho}
\ge L_{\rho, d}^{cl} \, {\tilde{u}_{1}}^{-2} \,
\left(z-\lambda_1\right)_{+}^{\rho+d/2}.
\end{equation}
\end{corollary}
We also have the following universal lower bound.
\begin{corollary}
\begin{equation}
\label{hermi-1} \sum_{k} \left(z-\lambda_k\right)_{+} \ge
\dfrac{2}{d+2} H_d^{-1} \lambda_1^{-d/2} \
\left(z-\lambda_1\right)_{+}^{1+d/2}.
\end{equation}
\end{corollary}
\par\noindent \textbf{Proof.} This corollary is evident using the
isoperimetric inequality of Chiti \cite{Chi3} \cite{Her},
\begin{equation}
\textup{ess sup} |u_1| \le
\left(\dfrac{\lambda_1}{\pi}\right)^{d/4}
\dfrac{2^{1-d/2}}{\Gamma(d/2)^{1/2} j_{d/2-1,1} J_{d/2}
(j_{d/2-1,1})}. \label{Chiti-ess-sup}
\end{equation}
With the way $H_d$ and $L_{0,d}^{cl}$ are defined in
(\ref{constant-riesz-low}) and (\ref{classical-constant-0}), we
prefer to put this inequality in the form
\begin{equation}
{\tilde{u}_{1}}^{2} \le H_d L_{0,d}^{cl}
{\lambda_1}^{d/2}.\label{chiti-myway}
\end{equation}
Substituting (\ref{chiti-myway}) into (\ref{laptev-1}) leads to (\ref{hermi-1}).
\hfill\qed

\par\noindent \textbf{Remarks.}
Theorem \ref{riesz-low-thm} can now be proved by either of 
two simple steps:  
\begin{itemize}
\item[(i)] 
Applying the Riesz
iteration to (\ref{hermi-1}) leads to (\ref{riesz-low}). 
\item[(ii)] 
Alternatively,
Theorem \ref{riesz-low-thm} follows from  Corollary
\ref{cor-univer} applying Chiti's inequality (\ref{chiti-myway}).
In \cite{Safarov} Safarov relied instead on a result of E.~B.~Davies
\cite{Dav1},
\begin{equation}
\Ess |u_1| \le e^{1/{8 \pi}}\, \la_1^{d/4}, \label{davies}
\end{equation}
to obtain a statement similar to Theorem \ref{riesz-low-thm}. The
use of Chiti's inequality (\ref{chiti-myway}),
which saturates when $\Omega$ is an
$d-$ball,
improves Safarov's
constant, particularly for large dimension $d$; see the discussion
in \cite{Her}.
\end{itemize}

As a corollary, we have the following lower bound for $Z(t)$

\begin{corollary} For $t\ge 0$
\begin{equation}
Z(t) \ge \dfrac{\Gamma(1+d/2)}{H_d} \, \dfrac{e^{-\lambda_1
t}}{\left(\lambda_1 \, t\right) ^{d/2}}. \label{partition-low}
\end{equation}
\end{corollary}
\par\noindent \textbf{Proof.} 
We reason as in the derivation of
Kac's ineq. (\ref{Kac-inequality}) from Berezin-Li-Yau
(\ref{laptev-weidl}).  Apply the Laplace transform to
(\ref{riesz-low}) to obtain
\[\dfrac{\Gamma(1+\rho)}{t^{1+\rho}} Z(t) \ge
H_d^{-1} \lambda_1^{-d/2} \ \dfrac{\Gamma(1+\rho)
\Gamma(1+d/2)}{\Gamma(1+ \rho+ d/2)} \ \dfrac{\Gamma(1+ \rho+
d/2)}{ t^{1+\rho+d/2}} \, e^{-\lambda_1 t}.\] 
Simplifying
results in the statement of the corollary. \hfill\qed

An immediate consequence of this corollary is the following
universal lower bound for the zeta function in terms of the
fundamental eigenvalue.

\begin{corollary} For $\rho>d/2$
\begin{equation}
\zeta_{spec}(\rho) \ge \dfrac{\Gamma(1+d/2)}{H_d} \,
\dfrac{\Gamma(\rho - d/2)}{\Gamma(\rho)} \,
\dfrac{1}{\lambda_1 ^{\rho}}. \label{zeta-low}
\end{equation}
\end{corollary}
\par\noindent \textbf{Proof.} This corollary is evident by applying the Mellin
transform
\[\zeta_{spec}(\rho) = \dfrac{1}{\Gamma(\rho)} \ \int_{0}^{\infty}
t^{\rho-1} Z(t) dt\] to the statement (\ref{partition-low}) and
observing that the definition of the $\Gamma$ function leads to
\[\dfrac{1}{\lambda^{\rho}} =\dfrac{1}{\Gamma(\rho)} \, \int_{0}^{\infty} e^{-\lambda
t} t^{\rho-1} dt.  \]
\hfill\qed 
\smallskip

\begin{figure}[htb!]
  \begin{center}
    \includegraphics[width=.95\textwidth, bb=0 0 320 197]{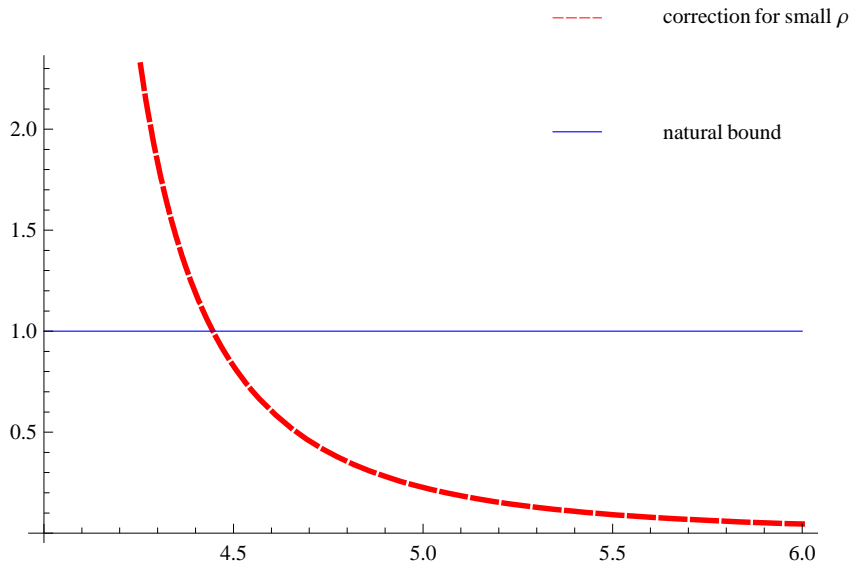}
    \caption{\textit{Universal Lower Bound Estimate for $\lambda_1^{\rho} \, \zeta_{spec}(\rho)$ from
    (\ref{zeta-low}) as a function of $\rho$, for $d=8$.}}
    \label{fig}
  \end{center}
\end{figure}

We also note that it is not hard to prove that there exists a
threshold value $\rho_0>d/2$ beyond which the estimate in
(\ref{zeta-low}) becomes weak (in comparison with dropping all the
terms in the definition of $\zeta_{spec}(\rho)$ except for
$1/\lambda_1^{\rho}$). This is illustrated in Fig. \ref{fig}.

Inequality \eqref{partition-low} lends itself to a generalization
in the spirit of Dolbeaut {\em et al.} \cite{DFLP}. We first adopt its setting. For a nonnegative function $f$ on $\rz_+$ such that
\[\int_0^{\infty} f(t) \left(1+ t^{-d/2}\right)
\frac{dt}{t}<\infty\] define \begin{equation} \label{weyl1}
F(s) :=\int_{0}^{\infty} e^{-s t}
f(t) \frac{dt}{t}\end{equation} and let 
\begin{equation}\label{weyl2}
G(s) := \mathcal{W}_{d/2} \{F(z)\}(s),
\end{equation}
where
\[\mathcal{W}_{\mu} \{F(z)\}(s):=\frac{1}{\Gamma(\mu)} \, \int_{s}^{\infty} F(z)
\left(z-s\right)^{\mu-1} dz\] denotes the Weyl transform of order
$\mu$ of the function $F(z)$. From the tables in \cite{bateman},
one notes that \[G(s)=\int_{0}^{\infty} \frac{e^{-s t}}{t^{d/2}}
f(t) \frac{dt}{t}.\] In fact, in analogy to what is shown in
\cite{DFLP}, \eqref{zeta-low} is a particular case of the
following. 

\begin{corollary} For $F(s)$ and $G(s)$ as defined above,
\begin{equation}
\sum_{j=1}^{\infty} F(\lambda_j) \ge \frac{\Gamma(1+d/2)}{H_d} \,
\lambda_1^{-d/2} G(\lambda_1). \label{general-lower}
\end{equation}
\end{corollary}
The proof of \eqref{general-lower} is immediate.
Scale \eqref{partition-low} by $f(t)/t$ then integrate from $0$ to
$\infty$. The counterpart to this inequality for Schr\"odinger operators 
has already been treated in \cite{DFLP}.

\par\noindent \textbf{Remarks.}

\begin{itemize}
\item[(i)] When $F(s)=s^{-\rho}$, $G(s)=\frac{\Gamma(\rho-\frac{d}{2})}{\Gamma(\rho)} \,
s^{d/2-\rho}$. Thus \eqref{zeta-low} is a particular case of \eqref{general-lower}. 

\item[(ii)] The choice $f(t)=a \, \delta(t-a)$, for $a>0$, leads to $F(s)=e^{-a s}$ and $G(s)=e^{-a s}/a^{d/2}$. One can then perceive that \eqref{partition-low} is a particular case of \eqref{general-lower} as well. 
Thus \eqref{partition-low} and \eqref{general-lower} are equivalent. 
\end{itemize}

\section{Remarks on the Work of A. Melas and Some Conjectures} \label{melas-work}

In  \cite{Melas} A. Melas proved the following inequality.
\begin{equation}
\sum_{i=1}^k \la_i \ge \dfrac{d}{d+2} \dfrac{4 \pi^2 k^{1+2/d}}
{\left(C_d |\Om|\right)^{2/d}} + M_d \, \frac{|\Omega|}{I(\Omega)}
\, k. \label{melas-inequality}
\end{equation}
Here $I(\Omega)$ is the ``second moment'' of $\Omega$, while $M_d$ is a constant that depends on the dimension $d$. Melas introduced the inequality as a correction to the Berezin-Li-Yau inequality \eqref{li-yau-sum}. 
\par\noindent


\par\noindent
Applying the Legendre transform $\Lambda\left[f\right](w) := \sup_z \left\{w z - f(z)\right\}$ (see \cite{LapWeidl2} \cite{LapWeidl1} \cite{hunder1} \cite{HarHer2})
to \eqref{melas-inequality}, one immediately obtains

\begin{equation}\label{new-melas1}
R_{\rho}(z) \le L_{\rho,d}^{cl} |\Omega|
\left(z-M_d
\frac{|\Omega|}{I(\Omega)}\right)_{+}^{\rho+\frac{d}{2}},
\end{equation}
for $\rho \ge 1$.  Applying the Laplace transform 
to \eqref{new-melas1} leads to the following correction
of Kac's inequality
\begin{equation}\label{new-melas2}
\sum_{i=1}^{\infty} \, e^{-\la_i t} \le \frac{|\Omega|}{\left(4
\pi t\right)^{d/2}} \, e^{-M_d \dfrac{|\Omega|}{I(\Omega)} \,
t}.
\end{equation}
Finally, applying the Weyl transform to \eqref{new-melas2}
leads to the following
\begin{equation} 
\zeta_{spec}(\rho) \le \frac{1}{(4 \pi)^{d/2}} \, \dfrac{
\Gamma(\rho-d/2)}{\Gamma(\rho)} \, |\Omega| \,  \left(M_d \dfrac{|\Omega|}{I(\Omega)}\right)^{\frac{d}{2}-\rho}. \label{zeta-high-melas}
\end{equation}
Furthermore,  reasoning as in Section \ref{lower-bounds},
these inequalities are particular cases of the following general theorem.
\begin{theorem}
For $F(s)$ and $G(s)$ as defined by \eqref{weyl1} and \eqref{weyl2}, one has
\begin{equation}
\sum_{j=1}^{\infty} F(\lambda_j) \le \frac{1}{(4 \pi)^{d/2}} \,
|\Omega| \, G\left(M_d \dfrac{|\Omega|}{I(\Omega)}\right). 
\label{general-upper-melas}
\end{equation}
\end{theorem}

We conjecture that a further improvement is possible, {\it viz.},
\begin{equation} \label{MelasConj}
\sum_{j=1}^{\infty} F(\lambda_j) \le \frac{1}{(4 \pi)^{d/2}} \,
|\Omega| \, G(|\Omega|^{-2/d})
\end{equation}
for the eigenvalues of the Dirichlet Laplacian, and that this is 
sharp.
In this case, 
$\frac{1}{|\Omega|^{2/d}}$ in \eqref{MelasConj} replaces 
$M_d \, \frac{|\Omega|}{I(\Omega)}$ in \eqref{general-upper-melas}.

Buttressing this conjecture is a related one for the spectral zeta function of the 
Dirichlet Laplacian:
\begin{conjecture}  \label{conj} For $\rho>d/2$,
\begin{equation} 
\zeta_{spec}(\rho) \le \dfrac{
\Gamma(\rho-d/2)}{\Gamma(\rho)} \ \dfrac{|\Omega|^{2
\rho/d}}{\left(4 \pi\right)^{d/2}}. \label{zeta-high-conj}
\end{equation}
\end{conjecture}
The conjectured universal constant 
\[C(\gamma)=\frac{1}{\left(4 \pi\right)^{d/2}} \, \frac{
\Gamma(\rho-d/2)}{\Gamma(\rho)} \ 
\]
appearing in this inequality is exactly that of the corresponding Schr\"odinger case in \cite{DFLP}. Statements \eqref{MelasConj} and \eqref{zeta-high-conj}
would be immediate consequences, using integral transforms,
of the following conjectured improvement to the Kac's inequality:
\begin{equation}\label{melas-conj2}
\sum_{i=1}^{\infty} \, e^{-\la_i t} \le \frac{|\Omega|}{\left(4
\pi t\right)^{d/2}} \, e^{-\dfrac{t}{|\Omega|^{2/d}}}.
\end{equation}

One might attempt to derive \eqref{zeta-high-conj} by emulating
\cite{LapWeidl2}, using a potential $V(x)$ equal to the characteristic function of
the complement of $\Omega$ multiplied by a coupling constant
tending to $+ \infty$,
but the constant that would appear on the right side of \eqref{zeta-high-conj}
is larger.  We point out that Conjecture \ref{conj} is consistent with
the Rayleigh-Faber-Krahn inequality
\[\lambda_1 \ge \dfrac{C_d^{2/d} \ j_{d/2-1,1}^2}{|\Omega|^{2/d}}\]
(as when one combines (\ref{zeta-low}) and
(\ref{zeta-high-conj})). Furthermore, as a result of (\ref{li-yau-sum}),
\begin{equation} \label{li-yau-eigen}
\lambda_k \ge \dfrac{d}{d+2} \dfrac{4 \pi^2 k^{2/d}} {\left(C_d
|\Om|\right)^{2/d}}, \notag
\end{equation}
and therefore
\begin{equation} \label{zeta-li-yau}
\zeta_{spec}(\rho) \le \left(\dfrac{d+2}{d}\right)^{\rho} \,
\dfrac{\zeta(2 \rho/d)}{\left(4 \pi^2\right)^{\rho}} \
\left(C_d \, |\Omega|\right)^{2\rho/d}.
\end{equation}

\begin{figure}[htb!]
  \begin{center}
  \caption{\textit{Upper Bound Estimate for $|\Omega|^{-2\rho/d} \, \zeta_{spec}(\rho)$ from
    (\ref{zeta-high-conj}), (\ref{zeta-li-yau}), and (\ref{zeta-polya}),
    as a function of $\rho$, for $d=2$.}}
    \label{fig2}
    \includegraphics[width=.85 \textwidth, bb=
0 0 320 197 ]{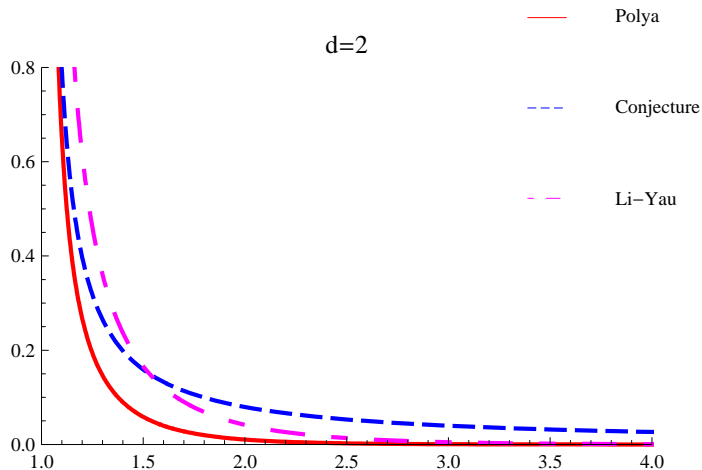}
  \end{center}
\end{figure}


If, as in the case of tiling domains, the P\'olya conjecture
\cite{PolyaConj}
\begin{equation}
\lambda_k \ge \dfrac{4 \pi^2 k^{2/d}} {\left(C_d
|\Om|\right)^{2/d}}. \label{polya-eigen} \notag
\end{equation}
is true, then
\begin{equation} \label{zeta-polya}
\zeta_{spec}(\rho) \le \dfrac{\zeta(2 \rho/d)}{\left(4
\pi^2\right)^{\rho}} \ \left(C_d \, |\Omega|\right)^{2\rho/d}.
\end{equation}


In both expressions above $\zeta$ denotes the usual expression for
the Euler zeta function, i.e.,
\[\zeta(\rho) = \sum_{k=1}^{\infty} \dfrac{1}{k^{\rho}}.\]

The bounds resulting from (\ref{zeta-high-conj}),
(\ref{zeta-li-yau}), and (\ref{zeta-polya}), for
$|\Omega|^{-2\rho/d} \, \zeta_{spec}(\rho)$ are plotted in
Fig. \ref{fig2}. It is clear that there is a threshold value
$\rho_0$ beyond which the conjectured bound
(\ref{zeta-high-conj}) cannot improve on Berezin-Li-Yau
(\ref{zeta-li-yau}). We expect that it should be 
possible to prove
\[\dfrac{\zeta(2 \rho/d)}{\left(4 \pi^2\right)^{\rho}} \, C_d^{2\rho/d}\le
\dfrac{1}{\left(4 \pi\right)^d} \,
\dfrac{\Gamma(\rho-d/2)}{\Gamma(\rho)} \le
\left(\dfrac{d+2}{d}\right)^{\rho} \, \dfrac{\zeta(2
\rho/d)}{\left(4 \pi^2\right)^{\rho}} \, C_d^{2\rho/d}.\]
Already Fig. \ref{fig2} gives credence to this statement and
Conjecture \ref{conj}.

\medskip

\textbf{Acknowledgements.} Our collaboration was initiated during the ``Low Eigenvalues of Laplace and Schr\"odinger Operators'' workshop held at the American Institute of Mathematics, Palo Alto (May 2006). The support of AIM is gratefully acknowledged. The second author would like to thank the Georgia Tech
School of Mathematics for their hospitality and support during his Fall 2006 visit. We also wish to thank Michael Loss and Joachim Stubbe for remarks and fruitful discussions.


\end{document}